\documentclass[reqno,centertags,12pt]{amsart}
\usepackage{amsmath,amsthm,amscd,amssymb,latexsym,verbatim}

\usepackage{graphicx,epsf,cite}

-\marginparwidth \oddsidemargin -4mm \evensidemargin \oddsidemargin

\newtheorem{theorem}{Theorem}[section]

\newtheorem{lemma}[theorem]{Lemma}
\newtheorem{corollary}[theorem]{Corollary}
\theoremstyle{definition}

\theoremstyle{remark}

\newcounter{smalllist}

\allowdisplaybreaks
\numberwithin{equation}{section}

\newcommand{\lb}{\label}

\newcommand{\beq}{\begin{equation}}
\newcommand{\eeq}{\end{equation}}

\newcommand{\bal}{\begin{align}}
\newcommand{\eal}{\end{align}}
\newcommand{\bals}{\begin{align*}}
\newcommand{\eals}{\end{align*}}

\newcommand{\bbR}{{\mathbb{R}}}

\newcommand{\bbT}{{\mathbb{T}}}

\newcommand{\til}{\tilde}

\begin{document}
\title[Rate of Merging of Vorticity Level Sets]
{On the Rate of Merging of Vorticity Level Sets \\ for the 2D Euler Equations}

\author{Andrej Zlato\v s}

\address{\noindent Department of Mathematics \\ University of
California San Diego \\ La Jolla, CA 92093 \newline Email: \tt
zlatos@ucsd.edu}


\begin{abstract} 
We show that two distinct level sets of the vorticity of a solution to the 2D Euler equations on a disc can approach each other along a curve at an arbitrarily large exponential rate.
\end{abstract}

\maketitle

\section{Introduction} \lb{S1}

In this note we study the question of how fast two distinct level sets of the vorticity of a solution to the Euler equations in two dimensions can approach each other.  We are here interested in the approach along a curve rather than just at a single point.  

The two-dimensional Euler equations model the motion of an incompressible inviscid fluid on a domain $D\subseteq\bbR^2$, and we will use here their vorticity formulation
\beq  \lb{1.1}
\omega_t+u\cdot\nabla\omega = 0
\eeq
on $D\times(0,\infty)$, with initial data
\beq  \lb{1.2}
\omega(\cdot,0) =\omega_0.
\eeq
We will consider the case when the vorticity $\omega=-\nabla \times u$ (which will be more convenient for us than the more standard $\omega=\nabla \times u$) is bounded, that is, $\omega_0\in L^\infty(D)$.
The customary no-flow boundary condition $u\cdot n=0$ on $\partial D\times(0,\infty)$, with $n$ the unit outer normal vector, then yields the Biot-Savart law
\[
u(x,t) =    - \int_{D}  \nabla^\perp G_D(x,y) \omega(y,t)dy
\]
for computing $u$ from $\omega$.  Here $\nabla^\perp=(-\partial_{x_2},\partial_{x_1})$ and $G_D$ is the Dirichlet Green's function for $D$ (i.e.,   $u=\nabla^\perp(- \Delta_D)^{-1}\omega$, with $\Delta_D$ the Dirichlet Laplacian on $D$).

It has been known since the works of H\" older \cite{Hol} and Wolibner \cite{Wol} that solutions to the Euler equations on smooth two-dimensional domains remain globally regular, and that $\| \nabla \omega(\cdot,t)\|_{L^\infty}$ cannot grow faster than double-exponentially  as $t\to\infty$ (although this bound seems to have first explicitly appeared in \cite{Yud}).  That is,  for each $\omega_0\in W^{1,\infty}(D)$ there is $C<\infty$ such that
\[
\|\nabla \omega(\cdot,t)\|_{L^\infty}\le C e^{e^{Ct}} \qquad\text{for each $t\ge 0$.}
\]
Whether the double-exponential rate of growth is attainable  had been a long-standing open problem.  The first examples of smooth solutions for which the vorticity gradient grows without bound as $t\to\infty$ were constructed by Yudovich \cite{Yud1, Yud2}.  Later Nadirashvili \cite{Nad} and Denisov \cite{Den} provided examples with at least linear and superlinear growth, respectively.  The period of relatively slow progress in this direction was ended by a striking recent result of Kiselev and \v Sver\' ak \cite{KS}.  Motivated by numerical simulations of Luo and Hou \cite{LuoHou2, LuoHou} which suggest blow-up for axisymmetric 3D Euler equations, they proved that solutions exhibiting double-exponential growth of the vorticity gradient indeed do exist in two dimensions.  This result, which was proved on a disc, was extended to general smooth two-dimensional domains with an axis of symmetry by Xu \cite{Xu}.   

The double-exponential growth in \cite{KS} is proved to occur on the boundary $\partial D$, whose presence is therefore crucial.   We note that the fastest growth currently known to occur on a domain without a boundary (i.e., on $\bbR^2$ or $\bbT^2$) is exponential growth of the vorticity gradient for solutions $\omega(\cdot,t)\in C^{1,\alpha}(\bbT^2)$ ($\alpha<1$), proved by the author in\cite{ZlaEulerexp} (see also \cite{LeiShi}).  Smooth solutions that grow super-linearly have been shown to exist as well, by Denisov \cite{Den}, who also constructed solutions exhibiting double-exponential growth rate on arbitrarily long but finite time intervals \cite{Den2} as well as patch solutions subject to a prescribed (regular) stirring for which the two patches approach each other double-exponentially in time \cite{Den3}.  
Finally, we note that  Kiselev and the author  proved that on domains whose boundaries are not everywhere smooth finite time blow-up can occur \cite{KZ}.

We will consider here \eqref{1.1} on a disc, as in \cite{KS}, although our results easily extend to general smooth two-dimensional  domains with  a symmetry axis via \cite{Xu}.  For convenience we will work with the unit disc $D:=B_1(e_2)$ centered at $e_2=(0,1)$, and we will denote its right/left halves by $D^\pm:=D\cap (\bbR^\pm\times\bbR)$.  Then we have
\[
G_D(x,y)=\frac 1{2\pi} \ln \frac{|x-y|}{|x-\bar y||y-e_2|}
\]
for $x\neq y\neq e_2$, with $\bar y := e_2+(y-e_2)|y-e_2|^{-2}$, as well as
\beq  \lb{1.3}
u(x,t) =    - \int_{D}  \left[ \frac{(x-y)^\perp}{|x-y|^2} - \frac{(x-\bar y)^\perp}{|x-\bar y|^2} \right] \omega(y,t)dy,
\eeq
where $(a_1,a_2)^\perp :=(-a_2,a_1)$.

The following is our main result.

\begin{theorem}\lb{T.1.1}
For $D=B_1(e_2)$ and each $A\ge 1$, there is $\delta>0$ and $\omega_0\in C^\infty (D)$ with $\|\omega_0\|_{L^\infty}=1$ such that the solution $\omega$ to \eqref{1.1}--\eqref{1.3} satisfies the following.  For any $t\ge 0$ we have
\[
\omega(0,\beta,t)=0 \qquad\text{for each $\beta\in(0,\delta)$,}
\]
and there is a function $\alpha_t:(0,\delta)\to(0,e^{-At})$ such that the set of those $\beta\in(0,\delta)$ for which
\[
\omega(\alpha_t(\beta),\beta,t)=1
\]
has measure at least $\delta- 2e^{-At}$.
\end{theorem}

{\it Remark.}   This result implies that $\|\omega(\cdot,t)\|_{W^{s,p}}$ also grows exponentially as $t\to\infty$ when $sp>1$.  On the other hand, the result in \cite{KS} yields double-exponential growth of these norms when $sp>2$ (as well as exponential growth for $(s,p)=(1,2)$).
\smallskip

The solutions that we will consider here are the ones from \cite{KS}, but we will track their dynamics in the neighborhood of the whole segment $\{0\}\times[0,\delta]$ rather than only near the origin.  A crucial extra ingredient in our argument will also be an {\it explicit} use of incompressibility of the flow $u$, and the corresponding measure-preserving property of its flow map (see \eqref{3.1}).  

We also note that if one were able to establish an additional estimate on these solutions, then one would obtain a super-exponential rate of merging of distinct level sets of $\omega$ (see Theorem~\ref{T.1.2} and Corollary \ref{C.1.3} below).  In the theoretically best possible case one could even prove a double-exponential rate of merging (see the discussion after Theorem \ref{T.1.2}), although it is not clear whether this case can occur.

In Section \ref{S2} we collect some estimates from \cite{KS} that we will use.  For the convenience of the reader and in order to provide more insight into the arguments that follow, we include the derivation of most of these, with the exception of the key Lemma \ref{L.2.1}.  We prove Theorem~\ref{T.1.1} in Section \ref{S3}, and the discussion of its extension when we are able to obtain additional estimates on $\omega$ appears in Section \ref{S4}.

\smallskip

{\bf Acknowledgements.}  
The author thanks Alexander Kiselev for useful discussions.  He also 
acknowledges partial support  by NSF grants DMS-1652284  and DMS-1656269.


\section{Some Estimates From \cite{KS}} \lb{S2}

To prove Theorem \ref{T.1.1}, we will employ some estimates from \cite{KS}, including  the following key lemma.  Just as in that paper, we will consider solutions $\omega$ that are odd in $x_1$ and non-negative on $D^+$.  That is, 
\beq\lb{2.1}
0\le \omega_0(x_1,x_2)=-\omega_0(-x_1,x_2)
\eeq
 for any $(x_1,x_2)\in D^+$.  Of course, then $\omega(\cdot,t)$ has the same properties for any $t\ge 0$ and $u_1(0,x_2,t)=0$ for any $(x_2,t)\in(0,2)\times(0,\infty)$.  Oddness of $\omega$ in $x_1$ also means that, with $\til x:=(-x_1,x_2)$ the reflection across the vertical axis, we have
\[
u(x,t) =    - \int_{D^+}  \left[ \frac{(x-y)^\perp}{|x-y|^2} - \frac{(x-\bar y)^\perp}{|x-\bar y|^2} - \frac{(x-\til y)^\perp}{|x-\til y|^2} + \frac{(x-\bar {\til y})^\perp}{|x-\bar {\til y}|^2} \right] \omega(y,t)dy.
\]
 We denote 
\begin{align*}
D_1^\gamma & :=\{x\in D^+ \,:\, x_1>\gamma x_2\},
\\ D_2^\gamma & :=\{x\in D^+ \,:\, x_2>\gamma x_1\}
\end{align*}
for $\gamma>0$, which are obtained by removing from $D^+$ sectors close to the $x_2$ and $x_1$ axes, respectively.  Finally, we  let 
\[
Q(x):=D^+\cap ([x_1,\infty)\times[x_2,\infty))
\]
for $x\in D^+$.

\begin{lemma}[\hskip -0.001mm \cite{KS}] \lb{L.2.1}
For any $\gamma>0$ there is $C_\gamma<\infty$ such that for each  $\omega_0\in L^\infty(D)$ that satisfies \eqref{2.1} we have
\beq\lb{2.2}
u_j(x,t)=   (-1)^j \left( \frac 4\pi\int_{Q(x)}\frac{y_1y_2}{|y|^4} \omega(y,t) dy  +  B_j(x,t)    \right) x_j 
\eeq
when $x\in D_j^\gamma$  ($j=1,2$), with $B_j$ satisfying
\beq\lb{2.3}
|B_j(x,t)|\le C_\gamma \|\omega_0\|_{L^\infty} .
\eeq
\end{lemma}

Note that when such $\omega_0$ is close to $\|\omega_0\|_{L^\infty}$ on all of $D^+$ except of  a small enough region (this property is then preserved by the evolution because $\omega$ is odd in $x_1$), then the first term in the parenthesis in \eqref{2.2} will dominate the second term for all $x$ close enough to the origin, regardless of where the small exceptional region is located.  Indeed,  if 
\beq\lb{2.4}
\left| \left\{ y\in D^+ \,:\, \omega_0(y)< \alpha \|\omega_0\|_{L^\infty} \right\} \right| \le \delta^2 \le \frac 1{16}
\eeq
for some $\alpha>0$ (this then also holds for $\omega(\cdot,t)$ and any $t\ge0$) and $x\in D^+\cap B_{2\delta}(0)$, then a simple analysis of the kernel $y_1y_2|y|^{-4}$  (which decreases radially and is maximized at $y_1=y_2$ for any fixed $|y|$) shows that
\[
\int_{Q(x)}\frac{y_1y_2}{|y|^4} \omega(y,t) dy \ge \alpha \|\omega_0\|_{L^\infty} \int_{[D^+\setminus B_{4\delta}(0)]\cap D_1^{1/2}\cap D_2^{1/2}} \frac{y_1y_2}{|y|^4} dy \ge \alpha\|\omega_0\|_{L^\infty}\int_{4\delta}^1 \int_{\pi/6}^{\pi/3} \frac{\sin 2\phi}{2r} d\phi dr
\]
for each $t\ge 0$.
Therefore, by Lemma \ref{L.2.1},
\beq\lb{2.5}
(-1)^j u_j(x,t)\ge \|\omega_0\|_{L^\infty}\left( \frac {\alpha |\ln 4\delta|}4 - C_\gamma   \right) x_j 
\eeq
for each $\gamma,\alpha,t>0$ and $j=1,2$ when \eqref{2.1} and \eqref{2.4} hold and $x\in D^+\cap B_{2\delta}(0)\cap D_j^\gamma$.  Of course, this estimate is only useful when $\delta>0$ is sufficiently small (depending on $\gamma,\alpha$).

We next consider the case $\|\omega_0\|_{L^\infty}=1$, for the sake of simplicity, although it is clear that the argument below works for any fixed $\|\omega_0\|_{L^\infty}>0$.  Let us take any $A\ge 1$, choose $\gamma=\tfrac 12$ and $\alpha=1$, and let 
\beq\lb{2.5a}
\delta:=\frac 14 e^{-4(A+C_{1/2})}.
\eeq
  If now $\omega_0$ satisfies  \eqref{2.1}, $\|\omega_0\|_{L^\infty}=1$,  and \eqref{2.4} with $\alpha=1$, then \eqref{2.5} yields for $j=1,2$,
\beq\lb{2.6}
(-1)^j u_j(x,t)\ge A x_j \qquad \text{for } (x,t)\in (D^+\cap B_{2\delta}(0)\cap D_j^{1/2})\times(0,\infty).
\eeq
 We pick such $\omega_0\in C^\infty(D)$, which also equals $1$ on the set $\{x\in D^+ \,:\, x_2\le x_1\in[\delta^{2}, \delta ] \}$.
 
Following \cite{KS}, for such $\omega_0$ and the corresponding solution $\omega$ and velocity $u$, let us denote 
 \begin{align*}
 \underline u_1(x_1,t) &:=\inf_{(x_1,x_2)\in D^+ \,\& \,x_2\le x_1} u_1(x_1,x_2,t), \\
 \overline u_1(x_1,t) &:=\sup_{(x_1,x_2)\in D^+ \,\&\, x_2\le x_1} u_1(x_1,x_2,t),
 \end{align*}
 and let $a,b$ be the solutions of 
 \begin{align} 
a'(t) &= \overline u_1(a(t),t), \qquad a(0)=\delta^{2}, \lb{2.21}
\\ b'(t) &= \underline u_1(b(t),t), \qquad b(0)=\delta.  \lb{2.22}
\end{align}
Note that $a,b$ are decreasing due to \eqref{2.6}, and they are positive on $(0,\infty)$ because  \eqref{2.2} yields the bound
\beq\lb{2.20}
|u_1(x,t)| \le \frac 4\pi \|\omega_0\|_{L^\infty} (\ln 2 - \ln |x|+C_\gamma) x_1
\eeq
for any $x\in D^+$ with $x_2\le x_1$.
Since \eqref{2.6} shows that 
\[
u_1(x,t)<0<u_2(x,t) \qquad\text{when $0<x_1=x_2\le\delta $ and $t>0$,}
\]
it follows that 
\beq\lb{2.5b}
\omega(\cdot,t)=1 \qquad \text{on $\{x\in D^+ \,:\, x_2\le x_1\in[a(t),b(t)] \}$}
\eeq
for any $t\ge 0$ such that $\inf_{s\in[0,t]}(b(s)-a(s))> 0$.  For any such $t$,
comparing \eqref{2.2} for points $(a(t),x_2)\in D^+$  with $x_2\le a(t)$ and for points $(b(t),y_2)\in D^+$  with $y_2\le b(t)$ (and also using the properties of $\omega(\cdot,t)$) yields 
\[
\frac d{dt} (\ln a(t) - \ln b(t)) \le -\frac 4\pi \int_{y\in D^+ \,\&\, a(t)<y_2<y_1\in [a(t),b(t)]} \frac {y_1y_2}{|y|^4} dy+ 
 \frac 4\pi  \int_{[b(t),1]\times[0,b(t)]} \frac {y_1y_2}{|y|^4} dy + 2C_{1/2}.
\]
The second integrand is  no more than $y_1^{-2}$ so the integral is bounded above by 1, and the first integral is bounded below by
\[
\int_{2a(t)}^{b(t)} \int_{\pi/6}^{\pi/4} \frac {\sin 2\phi}{2r} d\phi dr \ge \frac{\sqrt 3 \pi}{48}(\ln b(t) -\ln a(t) -\ln 2) \ge \frac{\pi}{32}(\ln b(t) -\ln a(t) -\ln 2),
\]
provided $b(t)\ge 2a(t)$.  This yields
\[
\frac d{dt} (\ln a(t) - \ln b(t)) \le \frac 18 (\ln a(t) - \ln b(t)) + C
\]
for such $t$, with $C:=\frac {\ln 2}8 + \frac 4\pi + 2C_{1/2}$.  Gronwall's lemma then shows 
\beq\lb{2.7}
\ln\frac{a(t)}{b(t)} \le \left( \ln\frac{a(0)}{b(0)}+ C \right) e^{t/8} 
\eeq
on any interval $[0,T]$ such that $\sup_{t\in[0,T]} \frac{a(t)}{b(t)}\le \frac 12$.  Since $\delta  <  e^{-1-C}$ due to \eqref{2.5a},  the parenthesis in \eqref{2.7} is less than $-1$ ($< \ln\frac 12$), and we thus obtain \eqref{2.7} for all $t\ge 0$.  We then have
\beq\lb{2.10}
a(t)\le \frac{a(t)}{b(t)}\le  (\delta e^C)^{e^{t/8}} \le e^{-e^{t/8}}
\eeq
for all $t\ge 0$.  In particular, $\|\nabla\omega(\cdot,t)\|_{L^\infty}$ grows double-exponentially in time (because  \eqref{2.5b} holds and $\omega(0,t)=0$ for all $t\ge 0$), which is the main result of \cite{KS}.

\section{Proof of Theorem \ref{T.1.1}} \lb{S3}

Consider now any $A\ge 1$, and let $\delta$ be from \eqref{2.5a} and $\omega_0$  as in the two sentences following \eqref{2.5a}.  Then $a,b$ from the previous section satisfy \eqref{2.5b} and \eqref{2.10}.
Fix any $T\ge 0$, let 
\[
V_T:=(0,a(T))^2\cap D^+,
\]
  and for $s\ge 0$ let
\[
U_{T}^s:=\Phi_{T}^s(V_T),
\]
where $\Phi$ is the flow map for $u$, given by
\[
\frac d{ds} \Phi_{T}^s(x) = u(\Phi_{T}^s(x), T+s), \qquad \Phi_{T}^0(x)=x\in D.
\]
That is, $U_T^s$ is the set to which the flow $u$ transports $V_T$ between times $T$ and $T+s$.  We also let
\[
\til U_T^s:=\{ \Phi_T^s (x) \,:\, x\in V_T \,\,\&\,\, \Phi_T^r (x)\in (0,\delta)^2 \text{ for each $r\in[0,s]$}  \}
\]
be the set of all points from $U_T^s$ which never left $(0,\delta)^2\cap D^+$ as they were transported by $u$ between times $T$ and $T+s$.
Since $u$ is divergence free, we have
\beq\lb{3.1}
|\til U_T^s| \le |U_T^s| =|V_T| < a(T)^2
\eeq
for each $T,s\ge 0$.

The estimate \eqref{2.6} applied at $\{x\in \overline{D^+}\,:\, 0\le x_{3-j}\le x_j=\delta\}$ shows that any point starting in $V_T$ at time $T$ and  leaving $(0,\delta)^2\cap D^+$ at a later time must do so across the line $[0,\delta]\times\{\delta\}$, as well as that it will also leave $[0,\delta]^2$ at the same time.  It therefore can lie in $\partial \til U_T^s$ at time $T+s$ only if that is the time of its first departure from $V_T$, in which case it also lies in $[0,\delta]\times\{\delta\}$.  It therefore follows that 
\[
\partial \til U_T^s \subseteq \Phi_T^s (\partial V_T) \cup ([0,\delta]\times\{\delta\})
\]
for each $T,s\ge 0$.
(Note that $\Phi_T^s$ extends continuously to $\partial D$ because both $\omega$ and $u$  extended continuously to $\partial D\times[0,\infty)$, with $u$ remaining log-Lipschitz because $\partial D$ is smooth.)

Let us now denote by $R_k$ ($k=1,2,3,4$) the four closed segments of $\partial V_T$ with endpoints at $(0,0)$, $(a(T),1-\sqrt{1-a(T)^2})$, $(a(T),a(T))$, and $(0,a(T))$, with $R_1$ being the segment between the first two of these points, and the other three segments labeled in counter-clockwise order.  Let us also fix 
\beq\lb{3.4}
s_T:=\frac 1A |\ln a(T)|.
\eeq
  Then \eqref{2.6} and $\delta\le \frac 14$ show that for each $x\in R_3$ there is $r< s_T$ such that $(\Phi_T^{r} (x))_2 = \frac 43\delta$.  In particular, 
\[
\Phi_T^{s_T} (R_3)\cap \partial \til U_T^{s_T} = \emptyset.
\]
Next, \eqref{2.6} together with $\Phi_T^{s_T}$ being continuous and satisfying  $\Phi_T^{s_T}(\overline{ D^+}\cap \overline{ D^-})=\overline{ D^+}\cap \overline{ D^-}$ and $\Phi_T^{s_T}(0)=0$ show that $ \{0\}\times [0,\delta] \subseteq \Phi_T^{s_T} (R_4)$ and in fact
\beq\lb{3.2}
 \{0\}\times [0,\delta] \subseteq \partial \til U_T^{s_T}
\eeq
From \eqref{2.6}, $\Phi_T^{s_T}(\partial D)=\partial D$ and $\Phi_T^{s_T}(0)=0$ we also have
\[
\Phi_T^{s} (R_1)\subseteq R_1
\]
for any $T,s\ge 0$.  From all this it follows that 
\beq\lb{3.3}
B:=\partial \til U_T^{s_T} \setminus \left[R_1\cup(\{0\}\times[0,\delta])\cup([0,\delta]\times\{\delta\})  \right] \subseteq \Phi_T^{s_T} (R_2).
\eeq
Moreover, \eqref{3.2} shows that for each $\beta\in (0,\delta)$ there is $\alpha_T^\beta>0$ such that 
\beq\lb{3.5}
(\alpha_T^\beta,\beta)\in \partial \til U_T^{s_T} \qquad\text{and}\qquad (0,\alpha_T^\beta)\times\{\beta\} \subseteq \til U_T^{s_T}.
\eeq
Then from \eqref{3.1} 
we conclude that
\beq\lb{3.7}
\left| \left\{\beta\in (0,\delta ) \,:\, \alpha_T^\beta < a(T) \right\}\right| \ge \delta -a(T) \qquad\text{for each $T\ge 0$.}
\eeq

Finally, the definition of $B$ and \eqref{3.5} show that 
\beq\lb{3.12}
(\alpha_T^\beta,\beta)\in B  \qquad\text{for each $\beta\in (1-\sqrt{1-a(T)^2},\delta)$,}
\eeq
 and then \eqref{3.3} with $1-\sqrt{1-a(T)^2}\le a(T)$
yield 
\beq\lb{3.6}
\omega(\alpha_T^\beta,\beta,T+s_T)=1  \qquad\text{for each $\beta\in \left( a(T),\delta \right)$ and $T\ge 0$.}
\eeq
Indeed, this holds because \eqref{1.1} is a transport equation and $\omega(x,T)=1$ for any $x\in R_2$ due to \eqref{2.5b}.
 Since $\omega(0,x_2,T+s_T)=0$ for any $x_2\in(0,2)$, the result follows once we notice that \eqref{2.10} and \eqref{3.4} yield
\beq\lb{3.8}
T+s_T\le 8\ln|\ln a(T)| + \frac 1A |\ln a(T)|.
\eeq

Indeed, since $a(\cdot)$ is continuous, for each 
\[
t\ge -\frac 1A\ln a(0) = \frac 2A |\ln\delta|
\]
 there is $T_t$ such that $T_t+s_{T_t}=t$.  If we also require $t\ge \frac 2A {e^{8A}}$, then either 
\beq\lb{3.8a}
|\ln a(T_t)| = As_{T_t} \ge e^{8A}
\eeq
or $T_t\ge \frac 1A e^{8A}$. In the latter case we use $e^{8A}\ge (8A)^2$ (due to $A\ge 1$) to conclude $T_t\ge 64A$, and then \eqref{3.8a} again follows from \eqref{2.10}.  Since $\ln x\le \frac x{8A}$ for $x\ge e^{8A}$, we now obtain from \eqref{3.8} that $t\le \frac 2A |\ln a(T_t)|$.  This yields $a(T_t)\le e^{-At/2}$, 
 which together with \eqref{3.7} 
implies 
  \beq\lb{3.9}
\left| \left\{\beta\in \left( 0,\delta \right) \,:\, \alpha_{T_t}^\beta \le e^{-At/2} \right\}\right| \ge \delta - e^{-At/2} \qquad\text{for each $t\ge \max\left\{ 2e^{8A},12(1+C_{1/2}) \right\}$,}
\eeq
where we also used $A\ge 1$ and 
\[
|\ln\delta|=4(A + C_{1/2})+\ln 4\le 6(A + C_{1/2}).
\]
On the other hand,  \eqref{3.6} and  $a(T_t)\le e^{-At/2}$
 now yield
\beq\lb{3.10}
\omega(\alpha_{T_t}^\beta,\beta,t)=1  \qquad\text{for each $\beta\in \left(e^{-At/2} ,\delta \right)$ and $t\ge \max\left\{ 2e^{8A}, 12(1+C_{1/2}) \right\}$.}
\eeq
Replacing $A$ by $2A$ now yields the result for all $t\ge \max\left\{ 2e^{16A}, 12(1+C_{1/2}) \right\}$.
To obtain it also for all $t\in [0,\max\left\{ 2e^{16A}, 12(1+C_{1/2}) \right\}]$, we only need to pick $\omega_0$ as above which is also equal to 1 on $D^+\cap ([c,1]\times\bbR)$ for a small enough $c>0$ (because $u$ is log-Lipschitz \cite{Hol, Wol}, so $|u_1(x,t)|\le C(|\log x_1|)x_1$ holds for some $C$ and all $(x,t)\in D^+\times(0,\infty)$).

\section{Towards a Faster Rate of Merging} \lb{S4}
 
 Notice that the double-exponential upper bound \eqref{2.10} on $a(t)$  is in fact not crucial for proving an exponential rate of approach of level sets of $\omega$ along a segment.  Indeed, if we only knew $a(T)\le e^{-cT}$ for  some $c>0$ and all $t\ge 0$, we would obtain 
 \[
 t=T_t+s_{T_t} \le \frac {A+c}{Ac} |\ln a(T_t)|.
 \]
   Then $a(T_t)\le e^{-Act/(A+c)}$, which would yield  \eqref{3.9} and \eqref{3.10} with $e^{-At/2}$ replaced by $e^{-Act/(A+c)}$ and a different lower bound on $t$. 
 
The limitation on obtaining a stronger result using \eqref{2.10} comes from our use of the bound \eqref{2.6} for $u_2$ on the time interval $[T,T+s_T]$, which dictates our choice of $s_T$ from \eqref{3.4}.  However, if we could gain more mileage from \eqref{2.2} by proving a faster growth of the second coordinate of $\Phi_T^s(x)$ for all $x\in R_3$, we might be able to improve Theorem \ref{T.1.1} further.  

This would be the case, for instance, if we could obtain a better {\it lower bound} on $b(t)$ than  \eqref{2.10}, such as $b(t)\ge e^{-3e^{t/8}}$ for large $t$.  More generally, let us assume that we have
\beq\lb{4.20}
a(t)\le e^{-ce^{t/C}} \qquad\text{and}\qquad b(t)\ge e^{-\tilde ce^{t/C}}  \qquad\text{for all large enough $t$},
\eeq
with some $C>0$ and $c>2\tilde c>0$.  Let us then, for the sake of simplicity, define $V_T$ in Section \ref{S3} with $e^{-ce^{T/C}}$ in place of $a(T)$, and change everything else up to the definition of $R_k$ accordingly (also, below only consider $T$ large enough so that \eqref{4.20} holds for all $t\ge T$).

Now for any $t\ge 0$ and any $a(t)\le x_1\le x_2\le b(t)$, we obtain
\[
u_2(x_1,x_2,t) \ge \left( \frac 4\pi\int_{x_2\le y_2\le y_1\le b(t)}\frac{y_1y_2}{|y|^4} dy  -C_{1/2}   \right) x_2
\]
from \eqref{2.2} and \eqref{2.5b}.  The integral can be estimated below by
\[
\int_{2x_2}^{b(t)} \int_{\pi/6}^{\pi/4} \frac {\sin 2\phi}{2r} d\phi dr \ge \frac{\sqrt 3 \pi}{48}(\ln b(t) -\ln 2x_2)
\] 
when $x_2\le \frac{b(t)}2$, so then
\beq\lb{4.2}
u_2(x_1,x_2,t) \ge \frac { -\ln x_2+\ln b(t) -\ln 2-8C_{1/2}}8 x_2.
\eeq
This is a better estimate than \eqref{2.6} when $x_2\le \frac {b(t)}2e^{-8(A+C_{1/2})}$, so it would be useful for $\Phi_T^s(x)$ when $x\in R_3$ as long as
\[
(\Phi_T^s(x))_2\le \frac{b(T+s)}{2e^{8(A+C_{1/2})}}.
\]
 
Let us therefore take $x\in R_3$ and consider $s\ge 0$ such that  $\Phi_T^s(x)$ did not yet exit $(0,\delta)^2$ --- thus, in particular, $(\Phi_T^s(x))_1\le (\Phi_T^s(x))_2$ --- and also such that with $\gamma:=\frac 12 e^{-8C_{1/2}}$ we have
\beq\lb{4.4}
(\Phi_T^s(x))_2\le  \gamma^2b(T+s)^2  \quad \left( \le \frac{b(T+s)}{2e^{8(A+C_{1/2})}} \text{ when $T\ge 4$ due to \eqref{2.5a} and \eqref{2.6}}\right).
\eeq
Then from \eqref{4.2} we obtain
\[
u_2(\Phi_T^s(x),T+s)\ge  \frac{|\ln (\Phi_T^s(x))_2|}{16}  (\Phi_T^s(x))_2.
\]
Since also $(\Phi_T^0(x))_2=x_2= e^{-ce^{T/C}}$, we obtain 
\beq\lb{4.9}
(\Phi_T^s(x))_2\ge e^{-ce^{(16T-Cs)/16C}} \qquad\text{until $(\Phi_T^s(x))_2=\gamma^2b(T+s)^2$.}
\eeq
The equality in \eqref{4.9} will be achieved at some  $s'\le \frac{16}C T$.
Then considering $s\ge s'$ and  applying \eqref{2.6} shows that $\Phi_T^s(x)$ exits $(0,\delta)^2$ at some 
\[
r<s'+\frac 2A |\ln \gamma b(T+s')| \le \frac{16}C T+ \frac 2A \left |\ln \gamma b \left(\frac{16+C}C T \right)\right|=:s_T.
\]
  The argument from the proof of Theorem \ref{T.1.1} then applies with this $s_T$, and we obtain \eqref{3.3}--\eqref{3.6} for all large $T$, with $e^{-ce^{T/C}}$ in place of $a(T)$.  Then for any large enough $t$, the equality $T_t+s_{T_t}=t$ from that proof becomes
\beq\lb{4.6}
\frac{16+C}C T_t+ \frac 2A \left |\ln \gamma b \left(\frac{16+C}C T_t \right)\right|=t.
\eeq
Since $b(t)\le e^{-At}$ by \eqref{2.6},  we have $\frac{16+C}CT_t\le \frac 1A |\ln \gamma b(\frac{16+C}CT_t)|$, so we know that $T_t$ satisfies
\beq\lb{4.3}
b\left(\frac{16+C}CT_t \right)   \le \gamma^{-1}e^{-At/3} =  2 e^{8C_{1/2}} e^{-At/3}.
\eeq
This finally yields  \eqref{3.9} and \eqref{3.10} with $e^{-At/2}$ replaced by $e^{-ce^{T_t/C}}$ and a different lower bound on $t$.  Note that this argument also requires 
\beq\lb{4.5}
a(t)\le \frac 14 e^{-16C_{1/2}} b(t)^2 \qquad\text{for all large enough $t$},
\eeq
 so that \eqref{4.4} holds at least for $s=0$ when $T=T_t$ (since obviously $\lim_{t\to\infty} T_t=\infty$), but this is guaranteed by $c>2\tilde c$ in \eqref{4.20}. (The power 2 in \eqref{4.4} and \eqref{4.5} could be replaced by any power $p>1$, at the expense of some constants above being different, so $c>\tilde c$ in fact suffices to obtain a result here.) 
 
 We have therefore proved the following result (notice that $b$ is decreasing, so replacing $\le$ by $=$ in \eqref{4.3} can only decrease $T_t$).
 
\begin{theorem}\lb{T.1.2}
Let $D=B_1(e_2)$, $D^+:=D\cap(\bbR^+\times\bbR)$, $C_{1/2}$ be from Lemma \ref{L.2.1},  $A\ge 1$, and  $\delta$ be from \eqref{2.5a}.
Consider any
$\omega_0\in C^\infty(D)$ satisfying  \eqref{2.1}, $\|\omega_0\|_{L^\infty}=1$, 
\[
\left| \left\{ x\in D^+ \,:\, \omega_0(x)< 1 \right\} \right| \le \delta^2,
\]
and equal to $1$ on the set $\{x\in D^+ \,:\, x_2\le x_1\in[\delta^{2}, \delta ] \}$, and let $\omega$ solve \eqref{1.1}--\eqref{1.3}.  If $a,b$ from \eqref{2.21}, \eqref{2.22} satisfy \eqref{4.20} with $C>0$ and $c>2\tilde c>0$, then the claims in Theorem \ref{T.1.1} hold for all large enough $t$, with $e^{-At}$ replaced by $e^{-ce^{T_t/C}}$ and  $T_t$ solving 
\[
b\left(\frac{16+C}CT_t \right)  =  2 e^{8C_{1/2}} e^{-At/3}.
\]
\end{theorem}

Let us now consider the theoretically best possible scenario.  For $\omega$ as above, the best lower bound on $b$ we could hope for is
\beq\lb{4.1}
b(t)\ge c' e^{-C't} \qquad\text{for all large enough $t$},
\eeq
with some $c'\in(0,1)$ and $C'\ge A$ (due to \eqref{2.6}).  If this is the case,
then \eqref{4.20} certainly holds and \eqref{4.3} yields
\[
T_t\ge  \frac {CA}{(50+3C)C'} t \qquad\text{for all large enough $t$}.
\]
That is,  \eqref{3.9} and \eqref{3.10} would hold with $e^{-At/2}$ replaced by $e^{-ce^{At/(50+3C)C'}}$ and a different lower bound on $t$.  Hence, the exponential lower bound \eqref{4.1} on $b(t)$ would yield a double-exponential rate of approach of two distinct level sets of $\omega$ along a segment.

Of course, it is not at all clear whether \eqref{4.1} can hold for some $\omega$ satisfying our hypotheses.  Nevertheless, even a much weaker lower bound may provide a stronger result than Theorem~\ref{T.1.1}.  For instance, if we can prove that
\beq\lb{4.21}
a(t)\le e^{-ce^{t/C}} \qquad\text{and}\qquad b(t)\ge e^{-c'e^{t/C'}}  \qquad\text{for all large enough $t$},
\eeq
with $c,c'>0$ and $C<C'$, then Theorem \ref{T.1.2} yields a super-exponential rate of approach whenever $C'>16+C$ (we obtain \eqref{3.9} and \eqref{3.10} with $e^{-At/2}$ replaced by $e^{-c(At/4C')^{C'/(16+C)}}$).

However, we can do even better in this case because if we only have a double-exponential lower bound on $b(t)$, then the above argument may be further optimized by this time choosing $s'$ so that $(\Phi_T^{s'}(x))_2=\gamma^2 e^{-2c'e^{(T+s')/C'}}$.
We then obtain $s'\le \frac{16(C'-C)}{(16+C')C}T+o(1)$
(with $o(1)=o(T^0)$ as $T\to\infty$). Therefore 
\[
T+s'\le \frac{(16+C)C'}{(16+C')C}T+o(1),
\]
so  we can replace \eqref{4.6} by
\[
\frac{(16+C)C'}{(16+C')C}T_t+o(1)+\frac 2A \left|\ln \gamma e^{-c'e^{[(16+C)C'T_t/(16+C')C+o(1)]/C'}}\right|=t.
\]
This yields
\[
T_t= \frac{(16+C')C}{16+C} \ln \frac{At}{2c'} + o(1),
\]
and hence leads to \eqref{3.9} and \eqref{3.10} with $e^{-At/2}$ replaced by $e^{-\kappa t^{(16+C')/(16+C)}}$ for some $\kappa >0$, and with a different lower bound on $t$.  So we again have a super-exponential rate of approach of distinct level sets of $\omega$ along a segment, but this time whenever $C<C'$:


\begin{corollary}\lb{C.1.3}
Let $D,A,\delta,\omega$ be as in Theorem \ref{T.1.2}.  If $a,b$ from \eqref{2.21}, \eqref{2.22} satisfy \eqref{4.21} with $c,c'>0$ and $C<C'$, then the claims in Theorem \ref{T.1.1} hold for all large enough $t$, with $e^{-At}$ replaced by $e^{-\kappa t^{(16+C')/(16+C)}}$ for some $\kappa >0$.
\end{corollary}





\end{document}